\documentclass[12pt]{amsart}
\usepackage{amsmath,amssymb,amsthm,amsfonts}
\usepackage{latexsym}
\usepackage{graphicx,color}
\def\vv<#1>{\langle#1\rangle}

\def\XXint#1#2{\setbox0=\hbox{$#1{#2}{\int}$}{#2}\kern-.5\wd0 }

\def\XXint#1#2#3{{\setbox0=\hbox{$#1{#2#3}{\int}$}
     \vcenter{\hbox{$#2#3$}}\kern-.5\wd0}}

\def\vv<#1>{\langle#1\rangle}

\newtheorem{theorem}{Theorem}[section]

\theoremstyle{definition}
\newtheorem{definition}{Definition}[section]
\theoremstyle{remark}

\numberwithin{equation}{section}

\begin{document}
\title{A geometric proof of the classification of complex vector cross product}

\author{Chengjie Yu}
\address{Department of Mathematics, Shantou University, Shantou,Guangdong, P.R.China}
\email{cjyu@stu.edu.cn}
\date{Oct. 2008}
\begin{abstract}
In this article, we give a geometric proof of the classification of complex vector cross product in Lee-Leung \cite{Lee-Leung}.
\end{abstract}
\maketitle
The following definition is a point-wise version of the definition of complex vector cross product on
complex Hermitian manifolds in Lee-Leung \cite{Lee-Leung}.
\begin{definition}
Let $V$ be complex vector space of dimension $n$ with a Hermitian metric and $r$ be an integer between $1$ and $n-1$. A $(r+1)$-form
$\phi\in \wedge^{r+1}V^*$ is called a complex $r$-fold vector cross product on $V$ if
\begin{equation}\label{eqn-cross-product}
\|\iota_{e_1\wedge e_2\wedge\cdots\wedge e_{r}}\phi\|^2=1
\end{equation}
for any $r$ unit vectors $e_1,e_2,\cdots,e_r$ that are orthogonal to each other.
\end{definition}

In  \cite{Lee-Leung}, Lee-Leung proved the following result.
\begin{theorem}
On a complex vector space $V$ of dimension $n$, there is a $r$-fold complex vector cross product if and
only if
\begin{enumerate}
\item $r=1$ and $n$ is an even number, or
\item $r=n-1$ and $n$ is an arbitrary natural number.
\end{enumerate}
\end{theorem}

In this article, we give a geometric proof of this theorem by using Pl\"uker embedding and dimension counting of projective
sub-varieties.
\begin{proof} The "if" part is clear (Ref. Lee-Leung \cite{Lee-Leung}). We only prove the "only if" part.

Let $\phi\in \wedge^{r+1}V^*$ be a $r$-fold complex vector product. Then, it induces a complex linear map
$f:\wedge^rV\to V^*$ by sending $u_1\wedge u_2\wedge\cdots\wedge u_r$ to $\iota_{u_1\wedge u_2\wedge\cdots\wedge u_r}\phi$.

By the identity (\ref{eqn-cross-product}), we know that for any $u_1,u_2,\cdots, u_r$ that are linearly independent,
\begin{equation*}
u_1\wedge u_2\wedge\cdots \wedge u_r\not\in \ker f.
\end{equation*}
By projectification and Pl\"uker embedding, all the $r$-vectors that can be separated as the form $u_1\wedge u_2\wedge\cdots \wedge u_r$
forms the Grassmannian manifold $Gr(r,V)\subset \mathbb{P}(\wedge^rV)$. The discussion above implies that that
\begin{equation*}
Gr(r,V)\cap \mathbb{P}(\ker f)=\emptyset.
\end{equation*}
Hence (Ref. Shafarevich \cite{Shafarevich})
\begin{equation*}
\dim Gr(r,V)+\dim \mathbb{P}(\ker f)<\dim \mathbb{P}(\wedge^r V).
\end{equation*}
That is,
\begin{equation}\label{eqn-1}
\dim Gr(r,V)+\dim \ker f<\dim (\wedge^r V).
\end{equation}
Note that, by dimension theorem in linear algebra, we have
\begin{equation}\label{eqn-2}
\dim \ker f\geq \dim \wedge^r V-n.
\end{equation}
Combining equation (\ref{eqn-1}) and (\ref{eqn-2}), we have
\begin{equation*}
r(n-r)=\dim Gr(r,V)<n.
\end{equation*}
Therefore $r=1$ or $r=n-1$.

When $r=1$, $f:V\to V^*$ is an isometry satisfying that
\begin{equation*}
\vv<f(u),u>=\phi(u,u)=0.
\end{equation*}
Let $e_1,e_2,\cdots,e_n$ be an orthonormal basis of $V$ and $\omega_1,\omega_2,\cdots,\omega_n$ be its dual
basis. Then, under this two basis, the matrix $A$ of $f$ is skew symmetric and nonsingular. Hence, $n$ is an even number.

\end{proof}

\end{document}